\documentclass[11pt,a4paper]{article}
\usepackage{latexsym}
\usepackage{graphicx}
\usepackage{amsmath}
\usepackage{amssymb}
\usepackage{epsfig}

\makeatletter \@addtoreset{equation}{section}
 \makeatother

\begin{document}

\title{ \textbf{On $\gamma$-Regular-Open Sets and $\gamma$-Closed Spaces}}
\author{SABIR HUSSAIN\\
Department of Mathematics, Islamia University Bahawalpur,
Pakistan.\\ {\bf Present Address:} Department of Mathematics,
Yanbu University,\\ P. O. Box 31387, Yanbu Alsinaiyah, Saudi
Arabia.\\E. mail: sabiriub@yahoo.com.}
\date{}
\maketitle
\textbf{Abstract.} The purpose of this paper is to continue
studying the properties of $\gamma$-regular open sets introduced
and explored in [6]. The concept of $\gamma$-closed spaces have
also been defined and discussed.
\\

\section {Introduction}

   The concept of operation $\gamma$ was initiated by S. Kasahara [7]. He
   also introduced $\gamma$-closed graph of a function. Using this
   operation,
H. Ogata [8] introduced the concept of $\gamma$-open sets and
investigated the related topological properties of the associated
topology $\tau_{\gamma}$ and $\tau$. He further investigated
general operator approaches of  close graph of mappings.

      Further S. Hussain and B. Ahmad [1-6] continued studying the properties of
   $\gamma$-open(closed) sets and generalized many classical notions in
   their work. The purpose of this paper is to continue
studying the properties of $\gamma$-regular open sets introduced
and explored in [6]. The concept of $\gamma$-closed spaces have
also been defined and discussed.

    First, we recall some definitions and results used in this paper. Hereafter, we shall
write a space in place of  a topological space.

\section {Preliminaries}

Throughout the present paper, X denotes topological spaces.
\\
{\bf Definition [7].} An operation $\gamma$ : $\tau\rightarrow$
P(X) is a function from $\tau$ to the power set of X such that V
$\subseteq V^\gamma$ , for each V $\in\tau$, where $V^\gamma$
denotes the value of $\gamma$ at V. The operations defined by
$\gamma$(G) = G, $\gamma$(G) = cl(G) and $\gamma$(G) = intcl(G)
are examples of operation $\gamma$.\\
{\bf Definition [7].}  Let $A \subseteq X$. A point $x \in A$  is
said to be $\gamma$-interior point of A, if there exists an open
nbd N of x such that $ N^\gamma\subseteq$ A and we denote the set
of all such points by $int_\gamma$(A). Thus
\begin{center}
        $int_\gamma$ (A) = $\{ x \in A : x \in N \in \tau$ and $N^\gamma\subseteq A \} \subseteq A$.
\end{center}

    Note that A is $\gamma$-open [8] iff A = $int_\gamma$(A). A set A is
called $\gamma$- closed [1] iff X-A is $\gamma$-open.\\
{\bf Definition [1].}  A point $x \in X$  is called a
$\gamma$-closure point of A$\subseteq$ X, if $U^\gamma\cap A \neq
\phi$, for each open nbd U of x. The set of all $\gamma$-closure
points of A is called $\gamma$-closure of A and is denoted by
$cl_\gamma$(A). A subset A of X is called $\gamma$-closed, if
$cl_\gamma(A) \subseteq A$. Note that $cl_\gamma(A)$ is contained
in every $\gamma$-closed superset
of A.\\
{\bf Definition [7].} An operation $\gamma$ on $\tau$ is said be
regular, if for any open nbds U,V of x $\in$ X, there exists an
open nbd W of x such that $U^\gamma \cap V^\gamma\supseteq
W^\gamma$.\\
{\bf Definition [8].} An operation $\gamma$ on $\tau$ is said to
be open, if for any open nbd U of each $x \in X$, there exists
$\gamma$-open set B such that $x \in B$ and  $U^{\gamma} \supseteq
B$.

\section {$\gamma$-Regular-Open Sets}

{\bf Definition 3.1 [6].} A subset A of X is said to be
$\gamma$-regular-open (respt. $\gamma$-regular-closed), if $A =
int_{\gamma}(cl_{\gamma}(A))$ (respt. $A =
cl_{\gamma}(int_{\gamma}(A)))$.

  It is clear that $RO_{\gamma}(X, \tau) \subseteq \tau_{\gamma}
  \subseteq \tau$ [6].\\

     The following example shows that the converse of above inclusion is not
     true in general.\\
{\bf Example 3.2.} Let X= $\{a,b,c\}$, $\tau =\{\phi , X, \{a\},
\{b\},\{a,b\},\{a,c\}\}$. For $b\in X$, define an operation
$\gamma:\tau \rightarrow P(X)$ by
\begin{center}
$\gamma(A) =  \left\{
\begin{array}{ccc}
A, & \mbox{if $b \in A$}\\
cl(A), & \mbox{if $b \not\in A$}\\
\end{array} \right.$\\
\end{center}
Calculations shows that $\{a,b\},\{a,c\},\{b\}, X ,\phi$ are
$\gamma$-open sets and  $\{a,c\},\{b\}, X ,\phi$ are
$\gamma$-regular-open sets. Here set $\{a,b\}$ is $\gamma$-open
but not $\gamma$-regular-open.\\
\\
{\bf Definition 3.3[7].} A space X is called $\gamma$-extremally
disconnected, if for all $\gamma$-open subset U of X,
$cl_{\gamma}(U)$ is a $\gamma$-open subset of X.\\
\\
{\bf Proposition 3.4.} If A is a $\gamma$-clopen set in X, then A
is a $\gamma$-regular-open set. Moreover, if X is
$\gamma$-extremally
disconnected then the converse holds.\\
{\bf Proof.} If A is a $\gamma$-clopan set, then $A =
cl_{\gamma}(A)$ and $A = int_{\gamma}(A)$, and so we have $A =
int_{\gamma}(cl_{\gamma}(A))$. Hence A is $\gamma$-regular-open.

  Suppose that X is a $\gamma$-extremally disconnected space and A is a
  $\gamma$-regular-open set in X. Then A is $\gamma$-open and so
  $cl_{\gamma}(A)$ is a $\gamma$-open set. Hence $A =
  int_{\gamma}(cl_{\gamma}(A))= cl_{\gamma}(A)$ and hence A is $\gamma$-closed
  set. This completes the proof.\\

   The following example shows that space X  to be $\gamma$-extremally
   disconnected is necessary in the converse of above Proposition.\\
{\bf Example 3.5} Let X= $\{a,b,c\}$, $\tau =\{\phi , X, \{a\},
\{b\},\{a,b\}\}$. Define an operation $\gamma:\tau \rightarrow
P(X)$ by $\gamma (B) = int (cl(B))$. Clearly X is not
$\gamma$-extremally disconnected space. Calculations shows that
$\{a\},\{a,b\},\{b\}, X ,\phi$ are $\gamma$-open as well as
$\gamma$-regular-open sets.  Here $\{a\}$ is a
$\gamma$-regular-open set but
not $\gamma$-clopan set.\\
\\
{\bf Theorem 3.6.} Let $A \subseteq X$, then $(a) \Rightarrow (b) \Rightarrow (c)$, where :\\
(a) A is $\gamma$-clopan.\\
(b) $A = cl_{\gamma}(int_{\gamma}(A))$.\\
(c) $X - A$ is $\gamma$-regular-open.\\
{\bf Proof.} $(a) \Rightarrow (b)$. This is obvious.

  $(b) \Rightarrow (c)$. Let $A = cl_{\gamma}(int{\gamma}(A))$.
  Then $X - A = X - cl_{\gamma}(int{\gamma}(A))= int_{\gamma}(X - int_{\gamma}(A)) = int_{\gamma}(cl_{\gamma}(X - A))$,
and hence $X -A$ is $\gamma$-regular-open set. Hence the proof.\\

   Using Proposition 3.4, we have the following Theorem:\\
{\bf Theorem 3.7.} If X is a $\gamma$-extremally disconnected
space. Then $(a) \Rightarrow (b) \Rightarrow (c)$, where :\\
(a) $X - A$ is $\gamma$-regular-open.\\
(b) A is $\gamma$-regular-open.\\
(c) A is $\gamma$-clopan.\\
{\bf Proof.} $(a) \Rightarrow (b)$. Suppose X is
$\gamma$-extremally disconnected space. From Proposition 3.4. , $X
-A$ is a $\gamma$-open and $\gamma$-closed set, and hence A is a
$\gamma$-open and $\gamma$-closed set. Thus $A =
int_{\gamma}(cl_{\gamma}(A))$ implies A is $\gamma$-regular-open
set.

   $(b) \Rightarrow (c)$. This directory follows from Proposition
   3.4. This completes as required.\\

   Combining Theorems 3.6 and 3.7, we have the following:\\
{\bf Theorem 3.8.}  If X is a $\gamma$-extremally disconnected
space. Then the following statements are equivalent:\\
(a) A is $\gamma$-clopan.\\
(b) $A = cl_{\gamma}(int_{\gamma}(A))$.\\
(c) $X - A$ is $\gamma$-regular-open.\\
(d) A is $\gamma$-regular-open.\\
\\
{\bf Theorem 3.9.} Let $A \subseteq X$ and $\gamma$ be an open
operation. If  $cl_{\gamma}(A)$ is a $\gamma$-regular-open set.
Then A is a $\gamma$-open set in X. Moreover, if X is extremally
$\gamma$-disconnected then the converse holds.\\
{\bf Proof.} Suppose that $cl_{\gamma}(A)$ is a
$\gamma$-regular-open sets. Since $\gamma$ is open, we have\\
$A \subseteq cl_{\gamma}(A) \subseteq
int_{\gamma}(cl_{\gamma}(cl_{\gamma}(A))) =
int_{\gamma}(cl_{\gamma}(A))= int_{\gamma}(A)$. This implies that
A is $\gamma$-open set.

  Suppose that X is $\gamma$-extremally disconnected and A is
  $\gamma$-open set. Then $cl_{\gamma}(A)$ is a $\gamma$-open set, and
  hence $\gamma$-clopan set. Thus by Theorem 3.8, $cl_{\gamma}(A)$
  is a $\gamma$-regular-open set. This completes the proof.\\
\\
{\bf Corollary 3.10.} Let X be a $\gamma$-extremally disconnected
space. Then for each subset A of X, the set
$cl_{\gamma}(int_{\gamma}(A))$ is $\gamma$-regular-open sets.\\
\\
{\bf Definition 3.11.} A point $x \in X$ is said to be a
$\gamma$-$\theta$-cluster point of a subset A of X, if
$cl_{\gamma}(U) \cap A \neq \phi$ for every $\gamma$-open set U
containing x. The set of all $\gamma$-$\theta$-cluster points of A
is called the $\gamma$-$\theta$-closure of A and is denoted by
$\gamma cl_{\theta} (A)$.\\
\\
{\bf Definition 3.12.} A subset A of X is said to be
$\gamma$-$\theta$-closed, if $\gamma cl_{\theta}(A) =A$. The
complement of $\gamma$-$\theta$-closed set is called
$\gamma$-$\theta$-open sets. Clearly a $\gamma$-$\theta$-closed
($\gamma$-$\theta$-open) is
   $\gamma$-closed($\gamma$-open) set.\\
\\
{\bf Proposition 3.13.} Let A and B be subsets of a space X. Then
the following properties hold:\\
(1) If $A \subseteq B$, then $\gamma cl_{\theta}(A) \subseteq
\gamma cl_{\theta} (B)$.\\
(2) If $A_i$ is $\gamma$-$\theta$-closed in X, for each $i \in I$,
then $\bigcap_{i \in I} A_i$ is $\gamma$-$\theta$-closed in X.\\
{\bf Proof.} (1). This is obvious.\\
(2). Let $A_i$ be a $\gamma$-$\theta$-closed in X for each $i \in
I$. Then $ A_i = \gamma cl_{\theta}(A_i)$ for each $i \in I$. Thus
we have\\
$\gamma cl_{\theta}(\bigcap_{i \in I}A_i) \subseteq \ \bigcap_{i
\in I} \gamma cl_{\theta}(A_i) = \bigcap_{i \in I} A_i \subseteq
\gamma cl_{\theta} (\bigcap_{i \in I} A_i)$.\\
Therefore, we have $\gamma cl_{\theta}(\bigcap_{i \in I} A_i) =
\bigcap_{i \in I} A_i$ and hence $\bigcap_{i \in I} A_i$ is
$\gamma$-$\theta$-closed. Hence the proof.\\
\\
{\bf Theorem 3.14.} If $\gamma$ is an open operation. Then for any
subset A of $\gamma$-extremally disconnected space X, the
following hold:
\begin{center}
$\gamma cl_{\theta}(A) = \bigcap \{$V : $A \subseteq V$ and V is
$\gamma$-$\theta$-closed$\}$\\
= $\bigcap \{$V : $A \subseteq V$ and V is
$\gamma$-regular-open$\}$
\end{center}
{\bf Proof.} Let $x \notin \gamma cl_{\theta}(A)$. Then there is a
$\gamma$-open set V with $x \in V$ such that $cl_{\gamma}(V) \cap
A = \phi$. By Theorem 3.9, $X - cl_{\gamma}(V)$ is
$\gamma$-regular-open and hence $X - cl_{\gamma}(V)$ is a
$\gamma$-$\theta$-closed set containing A and $x \notin X - \gamma
cl_{\theta}(V)$. Thus we have $x \notin  \bigcap \{$V : $A
\subseteq V$ and V is $\gamma$-$\theta$-closed$\}$.

   Conversely, suppose that $x \notin  \bigcap \{$V : $A
\subseteq V$ and V is $\gamma$-$\theta$-closed$\}$. Then there
exists a $\gamma$-$\theta$-closed set V such that $A \subseteq V$
and $x \notin V$,  and so there exists a $\gamma$-open set U with
$x \in U$ such that $U \subseteq cl_{\gamma}(U) \subseteq X - V$.
Thus we have $cl_{\gamma}(U) \cap A \subseteq cl_{\gamma}(U) \cap
V = \phi$ implies $x \notin \gamma cl_{\theta}(A)$.\\
 The proof of
the second equation follows similarly. This completes
the proof.\\
\\
{\bf Theorem 3.15.} Let $\gamma$ be an open operation. If X is a
$\gamma$-extremally disconnected
space and $A \subseteq X$. Then the followings hold:\\
(a) $x \in \gamma cl_{\theta}(A)$ if and only if $V \cap A \neq
\phi$, for each $\gamma$-regular-open set V with $x \in V$.\\
(b) A is  $\gamma$-$\theta$-open if and only if for each $x \in A$
there exists a $\gamma$-regular-open set V with $ x \in V$ such
that $V \subseteq A$.\\
(c) A is a $\gamma$-regular-open set if and only if A is
$\gamma$-$\theta$-clopan.\\
{\bf Proof.} (a) and (b) follows directly from Theorems 3.8 and 3.9.\\
(c) Let A be a $\gamma$-regular-open set. Then A is a
$\gamma$-open set and so $A = cl_{\gamma}(A) = \gamma
cl_{\theta}(A)$ and hence A is $\gamma$-$\theta$-closed. Since $X
- A$ is a $\gamma$-regular-open set, by the argument above, $X -A$
is $\gamma$-$\theta$-closed and A is $\gamma$-$\theta$-open. The
converse is obvious. Hence the proof.\\

  It is obvious that $\gamma$-regular-open $\Rightarrow$
  $\gamma$-$\theta$-open $\Rightarrow$ $\gamma$-open. But the
  converses are not necessarily true as the following examples
  show.\\
{\bf Example 3.16.} Let X= $\{a,b,c\}$, $\tau =\{\phi , X, \{a\},
\{b\},\{a,b\},\{a,c\}\}$. For $b\in X$, define an operation
$\gamma:\tau \rightarrow P(X)$ by
\begin{center}
$\gamma(A) =  \left\{
\begin{array}{ccc}
A, & \mbox{if $b \in A$}\\
cl(A), & \mbox{if $b \not\in A$}\\
\end{array} \right.$\\
\end{center}
Calculations shows that $\{a,b\},\{a,c\},\{b\}, X ,\phi$ are
$\gamma$-open sets as well as  $\gamma$-$\theta$-open sets and
$\gamma$-regular-open sets are $\{a,c\},\{b\}, X ,\phi$. Then the
subset $\{a,b\}$ is $\gamma$-$\theta$-open but not
$\gamma$-regular-open.\\
\\
{\bf Example 3.17.} Let X= $\{a,b,c\}$, $\tau =\{\phi , X, \{a\},
\{b\},\{a,b\},\{a,c\}\}$ be a topology on X. For $b\in X$, define
an operation $\gamma:\tau \rightarrow P(X)$ by
\begin{center}
$\gamma(A) = A^{\gamma}= \left\{
\begin{array}{ccc}
cl(A), & \mbox{if $b \in A$}\\
A, & \mbox{if $b \not\in A$}\\
\end{array} \right.$\\
\end{center}
Calculations shows that  $\{\phi , X, \{a\}, \{a,c\}\}$ are
$\gamma$-open sets and  $\{\phi , X, \{a,c\}\}$ are
$\gamma$-$\theta$-open sets. The the subset $\{a\}$ is
$\gamma$-open but not $\gamma$-$\theta$-open.

\section {$\gamma$-Closed Spaces}

{\bf Definition 4.1.} A filterbase $\Gamma$ in X,
$\gamma$-R-converges to $x_0 \in X$, if for each
$\gamma$-regular-open set A with $x_0 \in A$, there exists $F \in
\Gamma$ such that $F \subseteq A$.\\
\\
{\bf Definition 4.2.} A filterbase $\Gamma$ in X
$\gamma$-R-accumulates to $x_0 \in X$, if for each
$\gamma$-regular-open set A with $x_0 \in A$ and each $F \in
\Gamma$, $F \cap A \neq \phi$.\\

  The following Theorems directly follow from the above definitions.\\
{\bf Theorem 4.3.} If a filterbase $\Gamma$ in X,
$\gamma$-R-converges to $x_0 \in X$, then $\Gamma$
$\gamma$-R-accumulates to $x_0$.\\
\\
{\bf Theorem 4.4.} If $\Gamma_1$ and $\Gamma_2$ are filterbases in
X such that $\Gamma_2$ subordinate to $\Gamma_1$ and $\Gamma_2$
$\gamma$-R-accumulates to $x_0$, then $\Gamma_1$
$\gamma$-R-accumulates to $x_0$.\\
\\
{ \bf Theorem 4.5.} If $\Gamma$ is a maximal filterbase in X, then
$\Gamma$  $\gamma$-R-accumulates to $x_0$ if and only if
$\Gamma$  $\gamma$-R-converges to $x_0$.\\
\\
{\bf  Definition 4.6.} A space X is said to be $\gamma$-closed, if
every cover $\{V_{\alpha}: \alpha \in I \}$ of X by $\gamma$-open
sets has a finite subset $I_0$ of I such that $X = \bigcup_{\alpha
\in I} cl_{\gamma} (V_{\alpha})$.\\
\\
{\bf Proposition 4.7.} If $\gamma$ is an open operation, Then the
following are
equivalent:\\
(1) X is $\gamma$-closed.\\
(2) For each family $\{A_{\alpha} : \alpha \in I\}$ of
$\gamma$-closed subsets of X such that $\bigcap_{\alpha \in I}
A_{\alpha} = \phi$, there exists a finite subset $I_0$ of I such
that $\bigcap_{\alpha \in I_{0}} int_{\gamma}(A_{\alpha})=
\phi$.\\
(3)  For each family $\{A_{\alpha} : \alpha \in I\}$ of
$\gamma$-closed subsets of X, if $\bigcap_{\alpha \in I_{0}}
int_{\gamma}(A_{\alpha}) \neq \phi$, for every finite subset $I_0$
of I, then $\bigcap_{\alpha \in I} A_{\alpha} \neq \phi$.\\
(4) Every filterbase $\Gamma$ in X $\gamma$-R-accumulates to $x_0
\in X$.\\
(5) Every maximal filterbase $\Gamma$ in X $\gamma$-R-converges to
$x_0 \in X$.\\
{\bf Proof.} $(2) \Leftrightarrow(3)$. This is obvious.

     $(2)\Rightarrow (1)$. Let $\{ A_{\alpha}: \alpha \in I \}$ be a family
of $\gamma$-open subsets of X such hat $X = \bigcup_{\alpha \in I}
A_{\alpha}$. Then each $X - A_{\alpha}$ is a $\gamma$-closed
subset of X and $\bigcap_{\alpha \in I} (X - A_{\alpha}) = \phi$,
and so there exists a finite subset $I_0$ of $I$ such that
$\bigcap_{\alpha \in I_{0}} int_{\gamma}(X - A_{\alpha}) = \phi$,
and hence  $X = \bigcup_{\alpha \in I_{0}} (X - int_{\gamma}(X -
A_{\alpha}))= \bigcup_{\alpha \in I_{0}} cl_{\gamma}(A_{\alpha})$.
Therefore X is $\gamma$-closed, since $\gamma$ is open.

  $(4) \Rightarrow (2)$. Let $\{A_{\alpha} : \alpha \in I\}$ be a
family of $\gamma$-closed subsets of X such that $\bigcap_{\alpha
\in I} A_{\alpha} = \phi$. Suppose that for every finite subfamily
$\{A_{\alpha_i} : i = 1,2,...,n\}$, $\bigcap_{i=1}^{n}
int_{\gamma}(A_{\alpha_i}) \neq \phi$. Then $\bigcap_{i=1}^{n}
(A_{\alpha_i}) \neq \phi$ and  $\Gamma = \{\bigcap_{i=1}^{n}
A_{\alpha_i} : n \in N, \alpha_i \in I\}$ forms a filterbase in X.
By (4), $\Gamma$ $\gamma$-R-accumulates to some $x_0 \in X$. Thus
for every $\gamma$-open set A with $x_0 \in A$ and every $F \in
\Gamma$, $F \cap cl_{\gamma}(A) \neq \phi$. Since $\bigcap_{F \in
\Gamma} F = \phi$, there exists a $F \in \Gamma$ such that $x_0
\notin F$, and so there exists $\alpha_0 \in I$ such that $x_0
\notin A_{\alpha_0}$ and hence $x_0 \in X - A_{\alpha_0}$ and $X -
A_{\alpha_0}$ is a $\gamma$-open set. Thus $x_0 \notin
int_{\gamma}(A_{\alpha_0})$ and $x_0 \in X -
int_{\gamma}(A_{\alpha_0})$, and hence  $F_0 \cap (X -
int_{\gamma}(A_{\alpha_0})) = F_0 \cap cl_{\gamma}(X -
A_{\alpha_0}) = \phi$, which is a contradiction to our hypothesis.

  $(5) \Rightarrow (4)$. Let $\Gamma$ be filterbase in X. Then
  there exists a maximal filterbase $\xi$ in X such that $\xi$
  subordinate to $\Gamma$. Since $\xi$ $\gamma$-R-converges to
  $x_0$, so by Theorems 4.4
  and 4.5, $\Gamma$  $\gamma$-R-accumulate to $x_0$ .

  $(1)\Rightarrow (5)$. Suppose that $\Gamma = \{F_a :
  a \in I \}$ is a maximal filterbase in X which does not
  $\gamma$-R-converge to any point in X. From Theorem 4.5, $\Gamma$ does not $\gamma$-R-accumulates at any point in X.
  Thus for every $x \in
  X$, there exists a $\gamma$-open set $A_x$ containing x and $F_{a_x}
  \in \Gamma$ such that $F_{a_x} \cap cl_{\gamma}(A_x) = \phi$.
  Since $\{A_x : x \in X\}$ is  $\gamma$-open cover of X, there
  exists a finite subfamily $\{A_{x_i} : i = 1,2,...,n\}$ such
  that $X = \bigcup_{i=1}^{n} cl_{\gamma}(A_{x_i})$. Because
  $\Gamma$ is a filterbase in X, there exists $F_0 \in \Gamma$
  such that $F_0 \subseteq \bigcap_{i=1}^{n} F_{a_{x_i}}$, and
  hence $F_0 \cap cl_{\gamma}(A_{x_i})) = \phi$ for all $i =
  1,2,...,n$. Hence we have that, $\phi = F_0 \bigcap (\bigcup_{i=1}^{n} cl_{\gamma}(A_{x_i})) =
  F_0 \cap X$, and hence $F_0 = \phi$. This is a contradiction. Hence the proof.\\
\\
{\bf Definition 4.8.} A net $(x_i)_{i\in D}$ in a space X is said
to be $\gamma$-R-converges to $x \in X$, if for each $\gamma$-open
set U with $x \in U$, there exists $i_0$ such that $x_i \in
cl_{\gamma}(U)$ for all $i \geqslant i_0$, where D is a
directed set.\\
{\bf Definition 4.9.} A net $(x_i)_{i\in D}$ in a space X is said
to be $\gamma$-R-accumulates to $x \in X$, if for each
$\gamma$-open set U with $x \in U$ and each $i$, $x_i \in
cl_{\gamma}(U)$, where D is a
directed set.\\

 The  proofs of following Propositions are easy and thus are
 omitted:\\
{\bf Proposition 4.10.} Let $(x_i)_{i\in D}$ be a net in X. For
the filterbase $F((x_i)_{i\in D}) = \{ \{x_i : i \leq j \} : j \in
D
\}$ in X,\\
(1)$F((x_i)_{i\in D})$ $\gamma$-R-converges to $x$ if and only if
$(x_i)_{i\in D}$ $\gamma$-R-converges to $x$.\\
(2)$F((x_i)_{i\in D})$ $\gamma$-R-accumulates to $x$ if and only
if $(x_i)_{i\in D}$ $\gamma$-R-accumulates to $x$.\\
\\
{\bf Proposition 4.11.} Every filterbase F in X determines a net
$(x_i)_{i \in D}$ in X such that\\
(1) F $\gamma$-R-converges to $x$ if and only if $(x_i)_{i \in D}$
$\gamma$-R-converges to $x$.\\
(2) F $\gamma$-R-accumulates to $x$ if and only if $(x_i)_{i \in
D}$
$\gamma$-R-accumulates to $x$.\\

  From Propositions 4.11 and  4.12, filterbses and nets are
  equivalent in the sense of $\gamma$-R-converges and
  $\gamma$-R-accumulates. Thus we have the following Theorem:\\
{\bf Theorem 4.13.} For a space X, the following are equivalent:\\
(1) X is $\gamma$-closed.\\
(2) Each net $(x_i)_{i\in D}$ in X has a $\gamma$-R-accumulation
point.\\
(3) Each universal net in X $\gamma$-R-converges.


\begin{thebibliography}{99}
\bibitem{Ahmad}B. Ahmad and S. Hussain:\textit{ Properties of
$\gamma$-Operations on Topological Spaces, Aligarh Bull.
 Math.} 22(1) (2003), 45-51.
\bibitem{Ahmad}B. Ahmad and S. Hussain:\textit{ $\gamma$-Convergence in Topological Space ,
 Southeast Asian Bull. Math.,} 29(2005), 832-842.
 \bibitem{Ahmad}B. Ahmad and S. Hussain:\textit{ $\gamma^*$-Regular and $\gamma$-Normal Space ,
 Math. Today.,} 22(1)(2006), 37-44.
 \bibitem{Ahmad}B. Ahmad and S. Hussain:\textit{ On $\gamma$-s-Closed Subspaces ,
 Far East Jr. Math. Sci.,} 31(2)(2008), 279-291.
\bibitem{H}S. Hussain and B. Ahmad:\textit{ On Minimal $\gamma$-Open Sets,
Eur. J. Pure Appl. Maths.,} 2(3)(2009), 338-351.
\bibitem{H}S. Hussain and B. Ahmad:\textit{ On $\gamma$-s-Closed Spaces,
Sci. Magna Jr.,} 3(4)(2007), 89-93.
\bibitem{K}S. Kasahara:\textit{ Operation-Compact Spaces,
Math. Japon.,} 24(1979), 97-105.
\bibitem{H}H. Ogata:\textit{ Operations on Topological Spaces and Associated Topology,
Math. Japon.,} 36(1)(1991), 175-184.



\end{thebibliography}
\end{document}